\documentclass[12pt]{amsart}

\usepackage{amssymb}
\usepackage{amsthm}
\usepackage{amsmath}
\usepackage[all]{xy}
\usepackage[margin=1in]{geometry}

\usepackage{ulem}

\usepackage[usenames]{color}

\usepackage{tikz}
\usetikzlibrary{arrows,calc}

\usepackage{tikz-cd}

\usepackage[foot]{amsaddr}

\usepackage{hyperref}

\theoremstyle{plain}

\newtheorem{thm}     {Theorem}[section]

\newtheorem{definition}  [thm]{Definition}

\newtheorem{lemma}   [thm]{Lemma}

\newcommand{\C}{\mathbb C}
\newcommand{\D}{\mathbb D}
\newcommand{\N}{\mathbb N}

\newcommand{\R}{\mathbb R}

\begin{document}

\title{Properties and examples of Kobayashi hyperbolic Riemannian manifolds}

\author[H. Gaussier]{Herv\'e Gaussier}
\address{H. Gaussier: Univ. Grenoble Alpes, CNRS, IF, F-38000 Grenoble, France}
\email{herve.gaussier@univ-grenoble-alpes.fr}

\author[A. Sukhov]{Alexandre Sukhov$^1$}
\address{A. Sukhov: University  of Lille,   Laboratoire
Paul Painlev\'e,  Department of 
Mathematics, 59655 Villeneuve d'Ascq, Cedex, France, and
Institut of Mathematics with Computing Centre - Subdivision of the Ufa Research Centre of Russian
Academy of Sciences, 45008, Chernyshevsky Str. 112, Ufa, Russia.}
\email{sukhov@math.univ-lille1.fr}

\date{\today}
\subjclass[2010]{Primary: 32F45, 32Q45, 53A10, 53C15.}
\keywords{Kobayashi hyperbolic manifolds, Riemannian manifolds, conformal harmonic maps.}

\thanks{$^1\,$Partially suported by Labex CEMPI}

\begin{abstract}
We prove an analogue of the Brody lemma  in the framework of Riemannian manifolds. We also present new examples of Riemannian manifolds that are hyperbolic in the sense of Kobayashi.
\end{abstract}

\maketitle

\sloppy
\section*{Introduction}
The existence of the Kobayashi metric  in Riemannian manifolds was proved in \cite{Ga-Su}, generalizing results of B. Drinovec-Drnov\v sek - F.Forstneri\v c \cite{Dr-Fo}. The aim of this article is to study some properties of the Kobayashi metric. The classical Brody lemma provides a criterion for determining whether a compact complex manifold is Kobayashi hyperbolic, in terms of the existence of entire holomorphic curves. We present a partial analogue in Riemannian manifolds, based on a compactness result for harmonic maps. A Bloch principle was presented by M. Gromov \cite{Gr2}, dealing with harmonic maps in Riemannian manifolds. It states that sequences of harmonic discs have uniform bounded derivatives on a compact manifold as soon as there is no nonconstant harmonic map from $\R^2$ to $M$. Dealing with conformal harmonic maps, we provide a refined version of the Bloch principle, based on a Zalcman renormalization Lemma due to F. Berteloot \cite{Ber1}. This approach is very different from the approach in \cite{Gr2}.

Some examples of hyperbolic manifolds were proposed in \cite{Ga-Su2} where the notion of strict pseudoconvexity, first considered in \cite{Dr-Fo} for the case of the standard Euclidean metric in $\mathbb R^n$, was extended to  arbitrary Riemannian manifolds. We provide more examples of Kobayashi hyperbolic manifolds, considering Riemannian manifolds with negative sectional curvature.

The paper is organized as follows. In Section 1, we recall some results from \cite{Ga-Su, Ga-Su2}. In Section 2, we present the Brody Lemma. In Section 3, we prove that Riemannian manifolds with negatively pinched sectional curvature are Kobayashi hyperbolic.

\vspace{1mm}
This article is dedicated to the memory of Marek Jarnicki, for his remarkable contributions to the study of invariant metrics.

\section{Kobayashi hyperbolicity of Riemannian manifolds: preliminaries }

We denote by $\D$ the unit disc in $\R^2$ (we often identify $\R^2$ with the complex plane $\C$), by $ds^2$ the standard Riemannian metric on $\R^2$ and by $dm$ the standard Lebesgue measure on $\R^2$. For every $x \in \R^2$ and every $\lambda > 0$, we set $D(x,\lambda):=\{\zeta \in \C /\ |\zeta-x| < \lambda\}$. Throughout the paper, we consider a Riemannian manifold $(M,g)$, where $M$ is a smooth $C^{\infty}$ real manifold of real dimension $n \geq 2$ and $g$ is a smooth Riemannian metric of class $C^{\infty}$ on $M$. We denote by $dist_g$ the distance induced by $g$, defined as the infimum of the length of $C^1$ paths joining two points. For every $p \in M$ and every $r > 0$, we denote by $B_g(p,r)$ the ball $B_g(p,r) := \{q \in M /\ d_g(p,q) < r\}$.

%%%%%%%%%%%%%%%%%%%%%%%%%%%%%%%%%%%%%%%%%%%%

Let $(M,g)$ be a Riemannian manifold.  A map $u: \D \to M$ is called harmonic if it is a critical point of the energy integral

\begin{eqnarray*}
E(u) = \int_\D \vert du \vert^2dm.
\end{eqnarray*}

A smooth immersion $u : \D \to M$ is called conformal if  the pull-back $u^*g$ is a metric conformal to $ds^2$ i.e.,
there exists a smooth function $\phi$ such that $u^*g = e^\phi ds^2$ on $\D$.
A surface in $(M,g)$ is called minimal if its mean curvature induced by $g$ vanishes. A conformal  immersion (i.e., its image) is minimal if and only if it is harmonic (see \cite{Jo2}).

In the rest of the paper, we refer to conformal harmonic immersed disks as conformal harmonic immersions from $\D$ to $M$.

We denote $e_1:=(1,0) \in \R^2$. The following result claims the existence of a conformal harmonic disc with prescribed $1$-jet, see~\cite{Ga-Su2}, Lemma 2.1.

\begin{lemma}
\label{LemNW2}
 Let $p \in M$ and let $E$ be a 2-dimensional subspace of  $T_pM$. Then there exists a conformal harmonic immersion 
$u: \D \to M$ with $u(0) = p$ and such that the tangent space $T_p u(\D)$ coincides with $E$. Furthermore, this immersion depends smoothly on  $p$, $E$ and $g$.
\end{lemma}

Let $(M,g)$ be a Riemannian manifold. Using Lemma~\ref{LemNW2} we set, for a point $p \in M$ and a tangent vector   $\xi \in T_pM$
$$F_M(p,\xi) := \inf \frac{1}{r}$$
where $r$ runs over all positive real numbers for which there exists a conformal harmonic immersion  $u: \D \to M$ such that 
$u(0) = p$ and $du(0) \cdot e_1 = r \xi$. Then $F_M$ is well defined for every $(p,\xi)$ in the tangent bundle $TM$. We call $F_M$ the {\sl Kobayashi-Royden pseudometric} for the Riemannian manifold 
$(M,g)$.

The following properties of $F_M$ were proved in~\cite{Ga-Su}.

\begin{lemma}
\label{Lem1}
Let $f: (M,g) \to (N,h)$ be an isometric immersion   between two Riemannian manifolds i.e., $f$  satisfies $g = f^*h$. Then 
$$F_N(f(p),df(p) \cdot \xi) \le F_M(p,\xi).$$
In particular, if $M$ is a connected open subset of $N$, then 
$$F_N(p,\xi) \le F_M(p,\xi).$$
\end{lemma}

For $(p,\xi) \in TM$, we denote $\vert \xi \vert_g:=\sqrt{g(p)(\xi,\xi)}$, when no confusion is possible.
\begin{lemma}
\label{Lem2}
The function $F_M$ is non-negative, and for any real $a$ one has
$$F_M(p, a\xi) = \vert a \vert F_M(p,\xi).$$
If $K$ is a compact subset of $M$  then there is a constant $C_K > 0$ such that for $p \in K$ one has
$$F_M(p,\xi) \le C_K \vert \xi \vert_g.$$
\end{lemma}

\begin{thm}
\label{ThSemi}
The function $F_M$ is  upper semi-continuous on the tangent bundle $TM$.
\end{thm}

Denote by $P_{\D}$ the Poincar\'e metric, defined for $z \in \D$ and $v \in \C$ by

$$
P_{\D}(z,v) = \frac{\vert v \vert}{1- \vert z \vert^2}
$$

and by $\rho_{\D}$ the Poincar\'e distance on $\D$. Recall that for every $z,w \in \D$ :

$$
\rho_{\D}(z,w) = \frac{1}{2} \log \frac{1 + d(z,w)}{1- d(z,w)}
$$
where
$$
d(z,w) = \frac{\vert z - w \vert}{\vert 1 - z \overline w \vert}.
$$

 Let $p$ and $q$ be two points in the Riemannian manifold $(M,g)$. A Kobayashi chain from $p$ to $q$ is a finite
sequence of points $z_k$, $w_k$ in $\D$ and of conformal harmonic immersions $u_k: \D \to M$, $k = 1, \dots,m,$ such that $u_1(z_1) = p$, $u_m(w_m) = q$ and 
$u_k(w_k) = u_{k+1}(z_{k+1})$ for $k=1, \dots, m-1$. The Kobayashi pseudodistance from $p$ to $q$ is then defined by
\begin{eqnarray}
\label{KobDist1}
d_M(p,q) = \inf \sum_k \rho_{\D}(z_k,w_k)
\end{eqnarray}
where the infimum is taken over all Kobayashi chains from $p$ to $q$.

On another hand, we consider the pseudodistance defined as the integrated form of $F_M$:

\begin{eqnarray}
\label{KobDist2}
\overline d_M(p,q) = \inf \int_0^1 F_M(\gamma(t), \dot{\gamma}(t))dt
\end{eqnarray}
where the infimum is taken over all piecewise smooth paths $\gamma$ from $p$ to $q$ and where $\dot{\gamma}:=d\gamma/dt$. Notice that $\overline{d}_M$ is well defined by Theorem~\ref{ThSemi}. Then (see  \cite{Ga-Su}, Theorem 3.3)

\begin{thm}
\label{KobTh}
We have $d_M = \overline d_M$.
\end{thm}

\begin{definition}\label{hyp-def}
\begin{itemize}
\item[(i)] A Riemannian manifold $(M,g)$ is called {\rm Kobayashi hyperbolic} if the Kobayashi pseudodistance $d_M$ is a distance.
\item[$(ii)$] A Riemannian manifold $(M,g)$ is called {\rm hyperbolic} at a point $x \in M$ if there is a neighborhood $U$ of $x$ and a positive constant $c$ such that $F_M(y,\xi) \geq c \vert \xi\vert_g$ for every $y \in U$ and every $\xi \in T_yM$.

\item[$(iii)$] A Riemannian manifold $(M,g)$ is called {\rm tight} if the family of conformal harmonic immersions from $\D$ to $M$ is equicontinuous for the topology generated by $dist_g$.

\item[$(iv)$] A family $\mathcal F$ of mappings of a topological space $X$ into a topological space $Y$ is called  {\rm even} if, given $x \in X,\ y \in Y$ and a neighborhood $U$ of $y$, there is a neighborhood $V$ of $x$ and a neighborhood $W$ of $y$ such that for every $f \in F$, we have $f_{|V} \subset U$ whenever $f(x) \in W$.

\end{itemize}
\end{definition}

We recall the following characterization of hyperbolicity of Riemannian manifolds, see \cite{Ga-Su}
\begin{thm}
\label{HypTh}
Let $(M,g)$ be a Riemannian manifold. Then the following statements are equivalent:

\begin{itemize}
\item[$(i)$] the family $\mathcal C \mathcal H(\D,M)$ of conformal harmonic immersions from $\D$ to $M$ is equicontinuous with respect to the distance $dist_g$ i.e., $(M ,dist_g)$ is tight,

\item[$(ii)$] the family $\mathcal C \mathcal H(\D,M)$ is an even family,

\item[$(iii)$] $M$ is hyperbolic at every point,

\item[$(iv)$] $d_M$ is a distance i.e., $M$ is Kobayashi hyperbolic,

\item[$(v)$] the Kobayashi metric $d_M$ induces the usual topology of $M$.
\end{itemize}
\end{thm}

\section{Renormalization and the Bloch principle} 

\subsection{The Zalcman renormalization Lemma}

We follow the approach developed by F.Berteloot \cite{Ber1}. First, recall the following Schwarz property introduced in \cite{Ber1}.

\begin{definition}
Let $(J_\eta)_{\eta \in X}$ be a family of positive functions on a metric space $(X,d)$ such that $J_\eta(\eta) = 0$ for every $\eta \in X$. Let $0 < \alpha^- < \alpha^+ < 1$ 
and $ c > 0$ be some constants, and  $s: [0, \alpha^+[ \to \R^+$ be a function which is vanishing and continuous at $0$. We shall say that a map $f: \rho\D \to X$ satisfies 
the Schwarz property $S(J,\alpha^-,\alpha^+, c,s)$ if the following estimates occur for any $t_0, \kappa_0 \in \C$ such that $\vert t_0 \vert + \vert \kappa_0 \vert < \rho$:

\begin{eqnarray*}
\sup_{t \in \alpha^+\D} J_{f(t_0)}(f(t_0 + t\kappa_0)) \le c \implies J_{f(t_0)}(f(t_0 + t \kappa_0)) \le s(\vert t \vert)
\end{eqnarray*}
for each $t \in \alpha^-\D$.
\end{definition}

Here $\rho > 0$ and $\rho \D:=\{ t \in \mathbb C /\ |t| < \rho\}$.

The following version of the Zalcman renormalization lemma is obtained by F.Berteloot \cite{Ber1}.

\begin{lemma}
\label{Zalcman}
Let $(X,d)$ be a metric space and $(J_\eta)_{\eta \in X}$ be a family of positive functions on $X$ such that:
\begin{itemize}
\item[(i)] $J_\eta(\eta) =: 0$ for each $\eta \in X$,
\item[(ii)] $\lim_{\tau \to 0} \sup_{\eta \in K}(\sup_{d(\eta,\tilde \eta) \le \tau} J_\eta(\tilde \eta)) = 0$ for every compact subset $K$ of $X$.
\end{itemize}
Let $f_n: \overline\D \to X$  be a sequence of continuous maps such that:
\begin{itemize}
\item[(iii)]  every $f_n$ satisfies the Schwarz property $S(J,\alpha^-,\alpha^+,c,s)$ on $\D$,
\item[(iv)] there exist  $k \in ]0,1]$ and sequences $(\tilde t_n)_n$, $(\tilde\kappa_n)_n$ in $\C$ such that 
$\lim_n \tilde\kappa_n = 0$, $2\vert \tilde t_n \vert + \vert \tilde\kappa_n\vert < 1$  and $J_{f_n(\tilde t_n)} (f_n(\tilde t_n + \tilde \kappa_n) \ge k$ for 
every $n \in \N$.
\end{itemize}
Then there exists a sequence of affine contractions $r_n: t \mapsto t_n + \kappa_nt$ and a sequence $R_n > 0$ such that $\lim_n \kappa_n = 0$,
$\lim_n R_n = +\infty$, $\vert t_n \vert + \vert \kappa_n \vert < 1$, and
\begin{itemize}
\item[(a)]  $g_n := f_n \circ r_n$ is defined on $(R_n + \alpha^+)\D$ for every $n$,
\item[(b)]  $J_{g_n(0)}(g_n(1)) \ge k$ and $J_{g_n(t)}(g_n(t + u)) \le s(\vert u \vert)$ for every $n \in \N$, every $t \in R_n\D$ and every $u \in \alpha^-\D$.
\end{itemize}
Moreover, if $\lim_n \tilde t_n = 0$, then the contractions $r_n$ can be chosen such that $\lim_n r_n(0) = 0$.
\end{lemma}

\subsection{The Bloch principle and the Brody lemma}

Our goal now is to apply the Zalcman renormalization lemma to conformal harmonic maps. Let $(M,g)$ be a Riemannian manifold. Recall that an upper semi-continuous function $\rho : M \rightarrow [-\infty,+\infty)$ is minimal plurisubharmonic (MPSH) function if for every conformal harmonic immersed disc $u : \D \rightarrow M$, the composition $\rho \circ u$ is a subharmonic function on $\D$. The basic properties of MPSH functions are contained in \cite{Ga-Su2}.
The following lemma is proved in \cite{Ga-Su2}:

\begin{lemma}
\label{LemLog}
Suppose that the normal coordinates $x$ are chosen at the point $p = 0$.
Then the function $\phi(x) = \log \vert x \vert + A \vert x \vert$ is strictly MPSH in a neighborhood of the origin for a sufficiently large $A > 0$. In particular, for every $q$ sufficiently close to $p$ the function $ \phi_q(x) = \log \vert x - x(q) \vert + A \vert x - x(q) \vert$ is MPSH in the same neighborhood of the origin.
\end{lemma}
The size of a coordinate neighborhood and of a constant $A$ depend continuously on a point $p \in M$.

We assume now that $M$ is compact. Hence, we can choose $A$ independent of the point $p \in M$ and we can assume that each coordinate neighborhood contains a ball of radius equal to $3$ (up to scaling if necessary).  We fix such a neighborhood for every point $q \in M$.

Consider a smooth nondecreasing function $\psi$ on $\R_+$ satisfying $\psi(t) = t$ for $0 \le t \le 1/2$ and $\psi(t) = 1$ for $t \ge 3/4$. For each point $q \in M$ satisfying $\vert x(q) \vert < 2$, we define the function 
$\Psi_q =  \psi (\vert x - x(q) \vert^2) \exp(A\psi(\vert x - x(q) \vert)) $ on $D \cap U$, and 
$\Psi_q = 1$ on $D \setminus U$. Then the  function $\log \Psi_q (x) = \log \psi (\vert x - x(q) \vert^2) + A \psi(\vert x - x(q) \vert) + \lambda \rho$  is MPSH on $D \setminus \{ \vert x - x(q) \vert^2 \le 3/4 \}$. 

Obviously, $\Psi_q(q) = 0$ for every $q \in M$. Recall the following Schwarz lemma for subharmonic functions which is due to N.Sibony \cite{Si}.

\begin{lemma}
\label{Sibony}
Let $u$ be a function defined on $\D$, of class $C^2$ in a neighborhood of the origin. Suppose that $0 \le u \le 1$ , $u(0) = 0$ and $\log u$ is subharmonic on $\D$. Then
\begin{itemize}
\item[(a)] $u(z) \le \vert z \vert^2$ for every $z \in \D$ with equality at some point different from $0$ iff $u(z)$ is identically equal to $\vert z \vert^2$.
\item[(ii)] $\Delta u(0) \le 4$ with equality iff $u(z) = \vert z \vert^2$ for every $z \in D$.
\end{itemize}
\end{lemma}

Let us consider the family $(\Psi_q)$ (which we use instead of $(J_\eta)$ changing the notations). By construction, $(\Psi_q)$ satisfies Conditions (i) and (ii) of Lemma~\ref{Zalcman}. Moreover, it follows from (a) Lemma \ref{Sibony} that for the family $(\Psi_p)$ any conformal harmonic map $f: \D \to M$ satisfies the Schwarz property.

 We say that a smooth map $f: \D\to M$ is {\it weakly conformal} if $f^*g = \phi g_{st}$, with a smooth function $\phi \ge 0$. 

\begin{thm}
\label{Bloch}
Let $(M,g)$ be a compact Riemannian manifold. Assume that each weakly conformal harmonic map $f: \C \to M$ is constant. Then $(M,g)$ is Kobayashi hyperbolic.
\end{thm}
\proof Let $f_n: \D \to M$ be a sequence of conformal harmonic maps. It suffices to prove that this sequence is normal at the origin. Assume by contradiction that it is not. Since the family $(f_n)$ is not equicontinuous in any neighborhood of the origin,  Condition (iv) of Lemma \ref{Zalcman} holds. Hence, after renormalization, we obtain an equicontinuous family, which contains a  subsequence converging uniformly on compact subsets of $\C$ to a nonconstant map $f: \C \to M$. By the elliptic regularity, the convergence is in every  $C^k$ norm to a non-constant (in view of (b) Lemma \ref{Zalcman} ) weakly conformal harmonic map. This is a contradiction. \qed

\section{More examples of complete Kobayashi hyperbolic Riemannian manifolds}

As in Section 1, given two manifolds $M, N$, given a Riemannian metric $g$ on $M$ and a smooth map $f : N \rightarrow M$, we denote by $f^*g$ the $(0,2)$-tensor, pull-back of $g$ by $f$ and defined for every $p \in N$ and every $v,w \in T_pN$ by $f^*g(p)(v,w) = g(f(p))(df(p)\cdot v,df(p) \cdot w)$.
If $g$ is a Riemannian metric (resp. K\"ahler metric), we denote by $K(g)$ (resp. $H(g)$) the sectional (resp. holomorphic sectional) curvature of $g$. In order to avoid any confusion and since we will consider different Riemannian metrics on the same manifold, we will denote by $F_{(M,g)}$ the Kobayashi metric on $M$ with respect to the Riemannian metric $g$. First examples of complete Kobayashi hyperbolic Riemannian manifolds were given in \cite{Ga-Su2}. They include Riemannian manifolds $(M,g)$ with $K(g) = -c<0$.
The aim of the following is to provide more examples, studying the Kobayashi hyperbolicity of Riemannian manifolds $(M,g)$ with negative Riemannian sectional curvature, generalizing the case of hyperbolic Riemannian manifolds, and studying whether $g$ and its associated Kobayashi metric are bi-Lipschitz.

The computation of the  Kobayashi metric $F_{(M,g)}$ is based on the following version of the Schwarz Lemma for quasiconformal harmonic maps (see \cite{Go-Is} Theorem 2, p.565) :

{\bf Theorem \cite{Go-Is}.} {\sl Let $\D$ be endowed with the Poincar\'e metric $g_{\D}$ of constant curvature $-4$. Let $M$ be a Riemannian manifold with sectional curvature $K(g)$ bounded from above by a negative constant $-c$. If $u : \D \rightarrow M$ is a conformal harmonic mapping, then :
\begin{equation}\label{schwarz-eq}
u^*g \leq \frac{8}{c} g_{\D}.
\end{equation}
}
Here, for every $\zeta \in \D$, $\displaystyle g_{\D}(\zeta) = \frac{d\zeta \otimes d\overline{\zeta}}{(1-|\zeta|^2)^2}$.

\begin{definition}\label{quasi-def} Let $(M,h)$ be a Riemannian manifold.
For every $p \in M$, denote by  $\exp_p$ the exponential map defined in a neighborhood of $0$ in $T_pM$ and by $h_p$ the Riemannian metric $h_p:=(\exp_p)^*h$ (in particular, $h_p(0) = h(p)$). The Riemannian manifold $(M,h)$ has quasi-bounded geometry if there exists $r_0 > 0$ such that :

(a)  for every $p \in M$, $\exp_p$ is an immersion on $B_{h(p)}(0,r_0) \subset T_pM$, and

(b) for every integer $q \geq 0$, there exists $A_q$ such that for every $p \in M$, for every $|\mu| \leq q$ and for every $1 \leq i,j \leq n$,

$$
\sup_{x \in B_{h(p)}(0,r_0) \subset T_pM} \left|\frac{\partial^{|\mu|}((h_p-h(p))_{ij})}{\partial x^\mu}(x)\right| \leq A_q.
$$
\end{definition}

Then we have

\begin{thm}\label{hyp-is-hyp-thm}
Let $(M,g)$ be a Riemannian manifold with sectional curvature $K(g)$.

(i) If there exists $c > 0$ such that $K(g) \leq -c$, then $M$ is  Kobayashi hyperbolic.

(ii) If $(M,g)$ is complete with negatively pinched sectional curvature (i.e. there exists $c>0$ such that $-1/c \leq K(g) \leq -c$), then there exist a Riemannian metric $h$, bi-Lipschitz to $g$, and $C>0$ such that :
$$
\forall p \in M, \forall v \in T_pM, \frac{1}{C} \sqrt{h(p)(v,v)} \leq F_{(M,h)}(p,v) \leq C \sqrt{h(p)(v,v)}.
$$

(iii) If $(M,g)$ is complete with negatively pinched sectional curvature and with quasi-bounded geometry, then there exists $C>0$ such that :

$$
\forall p \in M, \forall v \in T_pM, \frac{1}{C} \sqrt{g(p)(v,v)} \leq F_{(M,g)}(p,v) \leq C \sqrt{g(p)(v,v)}.
$$
\end{thm}

Point (iii) of Theorem~\ref{hyp-is-hyp-thm} is a Riemannian analogue of Theorem 2 in \cite{Wu-Ya} in which the authors consider complete K\"ahler manifolds with negatively pinched holomorphic sectional curvature.
 
\begin{proof}
Point $(i)$ is a consequence of the Schwarz Lemma for Riemannian manifolds with sectional curvature bounded from above by a negative constant. Indeed, according to (\ref{schwarz-eq}), if $u : \D \rightarrow M$ is a conformal harmonic map, then for every $\zeta \in \D$:
$$
g(u(\zeta))(du(\zeta)\cdot e_1,du(\zeta)\cdot e_1) \leq \frac{8}{c} g_{\mathbb D}(\zeta)(e_1,e_1) = \frac{8}{c(1-|\zeta|^2)^2}.
$$
Hence, for every $p \in M$ and every $v \in T_pM$, if $u : \D \rightarrow M$ is a conformal harmonic map such that $u(0) = p, \ du(0) \cdot (e_1) = v/\alpha$ for some $\alpha > 0$, then
$$
\alpha \geq \sqrt{\frac{c}{8}} \sqrt{g(p)(v,v)},
$$
which implies
\begin{equation}\label{left-eq}
F_{(M,g)}(p,v) \geq \sqrt{\frac{c}{8}} \sqrt{g(p)(v,v)}.
\end{equation}

\vspace{2mm}
Let us consider now Point $(ii)$. Since the sectional curvature of $g$ is negatively pinched, we have (see Lemma 19 in \cite{Wu-Ya} or Theorem 1.2 in \cite{Sh} and Proposition in \cite{Ka})

\begin{lemma}\label{quasi-lem}
There exists a complete Riemannian metric $h$ on $M$ satisfying
$$
\left\{
\begin{array}{ll}
\displaystyle \frac{1}{C'} g \leq h \leq C'g, &\\
& \\
\displaystyle - \frac{1}{c'} \leq K(h) \leq -c' < 0,& \\
& \\
\displaystyle \sup_{x \in M} |\nabla_h^q R^h_{ijkl}| \leq C_q,& 
\end{array}
\right.$$
where $\nabla_h^q R^h_{ijkl}$ denotes the $q-th$ covariant derivatives of $R^h_{ijkl}$ with respect to $h$, and the positive constants $C' = C'(n)$, $c'=c'(n, c)$, $C_q = C_q(n, q, c)$ depend only on the parameters inside their parentheses.
\end{lemma}

It follows from Lemma~\ref{quasi-lem} that $(M,h)$ has quasi-bounded geometry, i.e., there exists $r_0 > 0$ such that :

(a)  for every $p \in M$, $\exp_p$ is an immersion on $B_{h(p)}(0,r_0) \subset T_pM$, and

(b) for every $p \in M$, let $h_p$ denote the Riemannian metric $h_p:=(\exp_p)^* h$, defined on $B_{h(p)}(0,r_0) \subset T_pM \simeq \mathbb R^n$ (in particular $h_p(0)=h(p)$). Then for every integer $q \geq 0$, there exists $A_q =A_q(n,q,c)> 0$ such that for every $p \in M$, for every $0 \leq |\mu| \leq q$ and for every $1 \leq i,j \leq n$,

\begin{equation}\label{bded-eq}
\sup_{x \in B_{h(p)}(0,r_0)} \left|\frac{\partial^{|\mu|}((h_p-h(p))_{ij})}{\partial x^\mu}(x)\right| \leq A_q.
\end{equation}

For every $t > 0$, let $\tilde{h}^t_p$ and $h^t_p$ be the Riemannian metrics defined for every $x \in B_{h(p)}(0,r_0/t)$ and every $v \in T_x(T_pM)\simeq T_pM$ by
$$
h^t_{p}(x)(v,v) := \frac{1}{t^2} \tilde{h}^t_{p}(x)(v,v) := h_p(t x)(v,v).
$$
Notice that $h^t_{p}$ is not isometric to $h_p$, but $\tilde{h}^t_{p}$ is isometric to $h_p$.
It follows from Condition (b) that for every $r > 0$ and every $k \in \mathbb N$, the metrics $h^t_{p}$ converge, when $t$ tends to zero, to the constant Riemannian metric $h_p(0)(=h(p))$ in $C^k(\overline{B_{h(p)}(0,r)})$ norm. Moreover, it also follows from Condition (b) that the convergence is uniform with respect to $p \in M$. From now on we set $r=2$.

For every $v, w \in T_pM \setminus\{0\}$ with $h(p)(v,w)=0$, let $u_{v,w}$ be the conformal harmonic disc (for $h(p)$) defined by
$$
\begin{array}{ccccc}
u_{v,w} & : & \D & \rightarrow & B_{h(p)}(0,1) \subset T_pM\\
 & & (x,y) & \mapsto & x \frac{v}{\sqrt{h(p)(v,v)}} + y \frac{w}{\sqrt{h(p)(w,w)}}\;.
 \end{array}
$$
Then $u_{v,w}$ satisfies $u_{v,w}(0)=0$ and $du_{v,w}(0) \cdot e_1 = v/\sqrt{h(p)(v,v)}$. Moreover, we have

\vspace{2mm}
\noindent{\bf Claim.}  {\sl There exist $0 < \varepsilon_0 << 1$ and $k_0 \geq 3$ satisfying the following : for every $p \in M$, for every smooth Riemannian metric $h'$ defined in a neighborhood of $\overline{B_{h(p)}(0,2)}$ and satisfying $\|h'-h(p)\|_{C^{k_0}(\overline{B_{h(p)}(0,2)})} < \varepsilon_0$, for every $v, w \in T_pM \setminus \{0\}$ with $h(p)(v,w) = 0$, there exists an immersed conformal harmonic disc $u'_{v,w} : \D \rightarrow B_{h(p)}(0,1)$, for the metric $h'$, such that

\begin{equation}\label{def-eq}
u'_{v,w}(0) = 0, \ T_0(u'_{v,w})(\D) = \mathbb R v \oplus \mathbb R w, \ \ du'_{v,w}(0) \cdot e_1 = \alpha \frac{v}{\sqrt{h(p)(v,v)}}, \ {\rm with} \ \alpha \geq \frac{1}{2}.
\end{equation}
}

The Claim follows from \cite{Ga-Su} Lemma~2.3 (i), where the immersion depends smoothly on the parameters $p, \zeta$ and on the metric, and from Condition~(\ref{bded-eq}) where the upper bound does not depend on $p$.

\vspace{2mm}
Let now $t_0 > 0$ be sufficiently small so that for every $p \in M$, $\|h^{t_0}_p - h(p)\|_{C^{k_0}\left(\overline{B_{h(p)}(0,2)}\right)} < \varepsilon_0$ (we also choose $t_0$ such that $t_0 < r_0/2$). Then for every $p \in M$, $v \in T_pM$ :

\begin{equation}\label{right-eq}
F_{(M,h)}(p,v) \underset{(a)}{\leq} F_{(\exp_p\left(B_{h(p)}(0,t_0),h\right)}(p,v) \underset{(b)}{\leq} F_{\left(B_{h(p)}(0,t_0),h_p\right)}(0,v) \underset{(c)}{\leq} \displaystyle \frac{1}{t_0}F_{\left(B_{h(p)}(0, 1),\tilde{h}^{t_0}_p\right)}(0,v)
\end{equation}

$$
 \underset{(d)}{\leq} \displaystyle \frac{1}{t_0} F_{\left(B_{h(p)}(0, 1),h^{t_0}_p\right)}(0,v) \underset{(e)}{\leq} \displaystyle \frac{2\sqrt{h(p)(v,v)}}{t_0}.
$$

Inequality (a) comes from the decreasing property of the Kobayashi-Royden metric under inclusion. Inequality (b) comes from the fact that by the definition of $h_p$ the map $\exp_p : (B_{h(p)}(0,t_0),h_p) \rightarrow (\exp_p\left( B_{h(p)}(0,t_0)\right),h)$ is a local isometry and satisfies $\exp_p(0) = p,\ d(\exp_p)(0) = Id_{T_pM}$. Inequality (c) comes from the fact that the map $\Lambda_0 : x \mapsto x/t_0$ is an isometry from $(B_{h(p)}(0,t_0),h_p)$ to $(B_{h(p)}(0,1),\tilde{h}^{t_0}_p)$. Inequality (d) comes from the condition $h^{t_0}_p=\frac{1}{t_0^2}\tilde{h}^{t_0}_p$ (in particular the set of stationary discs and hence of conformal harmonic immersed discs, for $\tilde{h}^{t_0}_p$ and $h^{t_0}_p$, coincide). Finally, Inequality (e) comes from Condition~(\ref{def-eq}).

\vspace{2mm}
Combining(\ref{left-eq}) and (\ref{right-eq}), we obtain finally :
$$
\forall p \in M, \ \forall v \in T_pM,\ \sqrt{\frac{c}{8}} \sqrt{h(p)(v,v)} \leq F_{(M,h)}(p,v) \leq \frac{2}{t_0}\sqrt{h(p)(v,v)}.
$$
This completes the proof of Point $(ii)$, setting $C=\max\left(2/t_0,\sqrt{c/8}\right)$.

\vspace{2mm}
Point $(iii)$ is a direct consequence of Point $(ii)$ with $h=g$. This concludes the proof of Theorem~\ref{hyp-is-hyp-thm}. \end{proof}


\begin{thebibliography}{CIT}

\bibitem{Ber1} {\sc Berteloot, F.}  {\sl Zalcman's lemma, Pinchuk’s rescaling method, and Catlin’s estimates revisited}, ArXiv : 2410.17633v1.

\bibitem{Dr-Fo}{\sc Drinovec-Drnov\v sek, B. ; Forstneri\v c, F.} {\sl Hyperbolic domains in real Euclidean spaces.} Pure and Applied Mathematics Quarterly {\bf 19} (2023), 2689-2735.

\bibitem{Fo-Ka}{\sc Forstneri\v c,F. ; Kalaj, D.} {\sl Schwarz-Pick lemma for harmonic maps which are conformal at a point.} ArXiv:2102.12403, to appear in Anal. PDE.

\bibitem{Ga-Su}{\sc Gaussier, H. ; Sukhov, A.} {\sl On the Kobayashi metrics on Riemannian manifolds.} ArXiv:2307.06154, to appear in Proc. Amer. Math. Soc.

\bibitem{Ga-Su2}{\sc Gaussier, H.; Sukhov, A.} {\sl Kobayashi hyperbolicity in Riemannian manifolds.} ArXiv:2407.15976.

\bibitem{Go-Is} {\sc Goldberg, S. I.; Ishihara, T.} {\sl Harmonic quasiconformal mappings of Riemannian manifolds.} Am. J. Math. {\bf 98} (1976), 225-240.

\bibitem{Gr2}{\sc Gromov, M.} {\sl Foliated Plateau problem, Part II: harmonic maps of foliations.} GAFA {\bf 1} (1991), 253-320.

\bibitem{Ha}{\sc Hawley, N.} {\sl Constant holomorphic curvature.} Canadian J. Math. {\bf 5} (1953), 53-56.

\bibitem{Ig}{\sc Igusa, J.} {\sl On the structure of a certain class of {K}aehler varieties.} Amer. J. Math. {\bf 76} (1954), 669-678.

\bibitem{Jo2}{\sc Jost, J.}  {\sl Riemannian geometry and geometric analysis.} Berlin: Springer-Verlag. xi, 401 p. Springer, 1995.

\bibitem{Ka}{\sc Kapovitch, V.} {\sl Curvature bounds via Ricci smoothing.} Ill. J. Math. {\bf 49} (2005),
259-263.

\bibitem{Ko} {\sc Kobayashi, S.} {\sl Hyperbolic complex spaces.} Grundlehren der Mathematischen Wissenschaften. 318. Berlin: Springer. xiii, 471 p. (1998). 

\bibitem{Sh}{\sc Shi, W.X.}{\sl Deforming the metric on complete Riemannian manifolds.} J. Differ. Geom. {\bf 30} (1989), 223-301.

\bibitem{Si} {\sc Sibony, N.}  {\sl A class of hyperbiolic manifolds}, in {\it Recent developments in several complex variables}, Ann. of Math. Study, vol. 100, Princeton Univ. Press, 
Princeton, N.J., 1981, pp. 357-372.

\bibitem{Wu-Ya}{\sc Wu, D.; Yau S.T.} {\sl  Invariant metrics on negatively pinched complete Kähler manifolds.} J. Am. Math. Soc. {\bf 33} (2020), 103-133.

\end{thebibliography}
\end{document}